\documentclass[11pt,a4paper,twoside]{amsart} %-*latex-*-
\usepackage{amsfonts}
\usepackage{amsmath}
\usepackage{amsthm}
\usepackage{amscd}
\usepackage{euscript}

\DeclareMathOperator{\HH}{H}

\DeclareMathOperator{\codim}{codim}

\newcommand{\Vop}{{ V_0 \otimes \dots \otimes V_p}}

\newcommand{\CC}{{\mathbb C}}

\newcommand{\Hh}{{\EuScript H}}

\newcommand{\OO}{{\EuScript O}}

\newcommand{\PP}{{\mathbb P}}

\newcommand{\fr}[1]{\stackrel{#1}{\longrightarrow}}

\newcommand{\ch}{\vee}
\newcommand{\Djs}{D_{j,\text{strong}}}

\swapnumbers
\newtheorem{thm}{Theorem}[section]
\newtheorem{lemma}[thm]{Lemma}
\newtheorem{cor}[thm]{Corollary}
\newtheorem{prop}[thm]{Proposition}
\theoremstyle{definition}
\newtheorem{defin}[thm]{Definition}

\theoremstyle{remark}
\newtheorem{remark}[thm]{Remark}
\newtheorem{exm}[thm]{Example}

\begin{document}
\title{Stabilizers for nondegenerate matrices of boundary format and Steiner bundles}
\author{Carla Dionisi}
\address{Carla Dionisi \\
 Dipartimento di Matematica Applicata
  ``G.~Sansone'' \\
  Via S.~Marta, 3\\
 I-50139, Firenze, Italy}
\email{dionisi@dma.unifi.it}


\begin{abstract}
  In this paper nondegenerate multidimensional matrices of boundary
  format in $V_0 \otimes \dots \otimes V_p$ are investigated by their
  link with Steiner vector bundles on product of projective spaces.
  For any nondegenerate matrix $A$ the stabilizer for the
  $SL(V_0)\times \dots \times SL(V_p)$-action, $Stab(A)$, is
  completely described. In particular we prove that 
  there exists an
  explicit action of $SL(2)$ on $V_0 \otimes \dots \otimes V_p$ such 
  that $Stab(A)^0\subseteq SL(2)$ and  the
  equality holds if and only if $A$ belongs to a unique $SL(V_0)\times
  \dots \times SL(V_p)$-orbit containing the identity matrices,
  according to \cite{ao}.
\end{abstract}
\maketitle

\markboth{{STABILIZERS FOR NONDEGENERATE MATRICES OF BOUNDARY FORMAT}}{CARLA DIONISI}
%%%%%%%%%%%%%%%%%%%%%%%%%%%%%%%%%%%%%%%%%%%%%%%%%%%%%%%%%
\section{Introduction}
Let $V_j$ be a complex vector space of dimension $k_j+1$ for
$j=0,\dots, p$ with $ k_0=\max_i \{k_i \}$. Gelfand, Kapranov and
Zelevinsky in \cite{gkz} proved that the dual variety of the Segre product
$\PP(V_0)\times\!\cdots\!\times\PP(V_p)$ is a hypersurface in
${(\PP^{(k_0+1)\dots(k_p+1)-1})}^\ch$ if and only if $k_0 \leq
\sum_{i=1}^p k_i$. The defining equation of this hypersurface is
called the {\it hyperdeterminant} of format $(k_0+1) \times \dots
\times (k_p +1)$ and is denoted by $Det$. Moreover the
hyperdeterminant is a homogeneous polynomial function on $V_0^\ch
\otimes \dots \otimes V_p^\ch$ so that the condition $Det A\neq0$ is
meaningful for a $(p+1)$-dimensional matrix $A\in \PP(V_0 \otimes
\dots \otimes V_p)$ of format $(k_0+1)\times \dots \times (k_p +1)$.
The hyperdeterminant is an invariant for the natural action of
$SL(V_0) \times \dots \times SL(V_p)$ on $\PP(V_0 \otimes \dots
\otimes V_p)$, and, in particular, if $Det A \neq 0$ then $A$ is
semistable for this action.

We denote by $Stab(A) \subset SL(V_0) \times \dots \times SL(V_p)$ the
stabilizer subgroup of $A$ and by ${Stab(A)}^0$ its
connected component containing the identity.
The stabilizer are well known for $p\leq1$ (in this case there is
always a dense orbit and the orbits are determined by the rank),
 so that in this paper we assume $p\geq2$. \\
It easy to check (see \cite{wz}, \cite{do}) that the degenerate matrices fill
an irreducible variety of codimension $k_0-\sum_{i=1}^pk_i+1$ and 
if $k_0<\sum_{i=1}^pk_i$ then all matrices  are degenerate.   
We will assume from now on that $A$ is of {\bf boundary format} i.e.,  that  $k_0=\sum_{i=1}^p k_i.$
 (A self-contained
approach to hyperdeterminant of boundary format matrices can be found
in \cite{do}).

\noindent  For multidimensional boundary
format matrices the classical definitions of triangulable,
diagonalizable and identity matrices can be easily reformulate in the
natural way  as follows
\begin{defin}
A $(p+1)$-dimensional matrix of boundary format $A \in V_0 \otimes
  \dots \otimes V_p$ is called 
\item {\bf triangulable} if $\forall j$ there exists a basis  $e_0^{(j)},
  \dots , e_{k_j}^{(j)} $ of $V_j$ such
  that \\ $ A=\sum a_{i_0,\dots , i_p}e_{i_0}^{(0)} \otimes \dots
\otimes e_{i_p}^{(p)}$ where $a_{i_o, \dots ,
    i_p}=0$ for $i_0 > \sum_{t=1}^p i_t $;
\item {\bf diagonalizable} if there exists a basis  $e_0^{(j)},
  \dots , e_{k_j}^{(j)} $ of $V_j$ such
  that \\ $ A=\sum a_{i_0,\dots , i_p}e_{i_0}^{(0)} \otimes \dots
\otimes e_{i_p}^{(p)}$ where $a_{i_o, \dots ,
    i_p}=0$ for $i_0 \neq \sum_{t=1}^p i_t$;
\item an {\bf identity} if there exists a basis  $e_0^{(j)},
  \dots , e_{k_j}^{(j)} $ of $V_j$ such
  that \\ $A=\sum a_{i_0,\dots , i_p}e_{i_0}^{(0)} \otimes \dots
\otimes e_{i_p}^{(p)}$ where 
\begin{equation*}
 a_{i_o, \dots , i_p}= \begin{cases}
0 & \text{for}\quad i_0 \neq \sum_{t=1}^p i_t \\
1  & \text{for} \quad i_0 = \sum_{t=1}^p i_t \\
\end{cases}
\end{equation*}
\end{defin}

Ancona and Ottaviani in \cite{ao}, considering the natural action 
of $SL(V_0) \times \dots
\times SL(V_p)$ on $\PP (V_0\otimes\ldots\otimes V_p)$, analyze these
properties from the point of view  of Mumford's Geometric Invariant Theory.

In the same aim, the main result of this paper is the following: 
\begin{thm}
\label{th:main}
Let $A\in \PP(V_0 \otimes \dots \otimes V_p)$ be a boundary format
matrix with $Det A \neq 0$.Then there exists a 2-dimensional vector
space $U$ such that $SL(U)$ acts over $V_i \simeq S^{k_i}U$ and
according to this action on $V_0 \otimes \dots \otimes V_p$ we have
$Stab (A)^0 \subseteq SL(U)$. Moreover the following cases are possible
$$Stab(A)^0\simeq \left\{
  \begin{matrix}
   0&\\\CC\\\CC^*&\\
SL(2)& \text{\ (this case occurs if and only if\ }A\hbox{\ is an
identity)}
  \end{matrix}\right.$$

\end{thm}
\begin{remark}
We emphasize that $SL(V_0) \times \ldots \times SL(V_p)$ is a "big"
group, so it is quite surprising that the stabilizer found lies always in the $3$-dimensional group $SL(U)$ without any dependence on $p$ and on dim$V_i$.

The maximal stabilizer is obtained by the "most symmetric" class of matrices 
corresponding to the identity matrices.
Under the identifications $V_i=S^{k_i}U$ the identity is given by the natural map
\begin{equation*}
  S^{k_1} U \otimes \ldots \otimes S^{k_p}U \to S^{k_0}U 
\end{equation*}
 which is defined under the assumption $k_0=\sum k_i$.
This explains again why the  condition of boundary format is so important.
\end{remark}
\noindent Ancona and Ottaviani in \cite{ao} prove theorem \ref{th:main} for
$p=2.$ We generalize their proof by using the correspondence
between nondegenerate boundary format matrices and vector bundles on a
product of projective spaces.  

\noindent Indeed, for any fixed $j \neq 0$, a
$(p+1)$-dimensional matrix $A \in V_0 \otimes \dots \otimes V_p$ of
format $(k_0+1)\times \dots \times (k_p +1)$ defines a sheaf morphism
${f_A}^{(j)}$ on the product $X=\PP^{k_1}\times \dots \times
\widehat{\PP^{k_j}}\times \dots \times \PP^{k_p}$
\begin{equation}
\label{steinerres}
 \OO_X \otimes V_0^{\ch} \fr{{f_A}^{(j)}}\OO_X(1, \dots, 1)\otimes V_j;
 \end{equation}
and it is easy to prove the following  
\begin{prop}
  \label{prop:det}({\cite{ao},\cite{cd}})
 If a matrix $A$ is of boundary format, then   $Det A \neq 0$ if and
only if  for all $j \neq 0$ the morphism ${f_A}^{(j)}$ is surjective (so ${{S}_A^\ch}^{(j)}=Ker {f_A}^{(j)}
$ is a  vector bundle of rank $k_0-k_j$). 
\end{prop}
In the particular case $p=2$ the (dual) vector bundle
${S_A}^{(1)}$ (or ${S_A}^{(2)}$) lives on the projective space $\PP^n
,\  n=k_2$ (or $ n=k_1$) and it is
a Steiner bundle as defined in \cite{dk} (this case has been
investigate in \cite{ao}). We shall refer to $S_A^{(j)}$ with the
name Steiner also for $p \geq 3$.

\noindent The main new technique introduced in this paper is the use of
jumping hyperplanes for bundles on the product of $(p-1)$ projective spaces.
For $p \geq 2$ there are two natural ways to introduce them;  by the above
correspondence, they translate  into two different conditions on the associated
matrix and that we call weak and strong (see definition \ref{def:wjumping} and
\ref{def:sjumping}). They coincide when $p=2$. 

Moreover, the loci of weak and strong jumping hyperplanes are
invariant for the action of $SL(V_0) \times \ldots \times SL(V_p)$ on
matrices.  By investigating these invariants we derive the proof of
theorem \ref{th:main} and also we obtain a characterization of a
particular class of bundles called Schwarzenberger bundles (see
\cite{s} for the original definition in the case $p=2$).
Schwarzenberger bundles correspond exactly to such matrices $A$ which
verify the equality $Stab(A)^0=SL(2)$ in theorem \ref{th:main}, called
identity matrices.

I would like to thank G.~Ottaviani  for
his invaluable guidance and the referee for useful suggestions to improve
this note.  

%%%%%%%%%%%%%%%%%%%%%%%%%%%%%%%%%%%%%%%%%%%%%%%%%%%%%%%%%
\section{Jumping hyperplanes and stabilizers}
%%%%%%%%%%%%%%%%%%%%%%%%%%%%%%%%%%%%%%%%%%%%%%%%%%%%%%%%%%%%5
Let  $p=2$ and $S:=S^1$ be the Steiner bundle on $\PP(V_2)$ defined
by a matrix $A\in V_0 \otimes V_1 \otimes V_2$ of boundary format, an hyperplane
$h \in \PP(V_2^\ch)$ is an unstable hyperplane of $S$ if
$h^0(S^\ch_{\vert h}) \neq 0$ (see \cite{ao}). By abuse of notations
we identify an hyperplane $h \in \PP(V_2^\ch)$ with any vector $h'\in
V_2$ such that $<h'>=h$.    

 In particular, $\HH^0(S^\ch(t))$ identifies to the
space of $(k_0+1) \times 1$-column  vectors $v$
with entries in $S^tV_2$ such that $Av=0$, and a
hyperplane $h$ is unstable  for $S$ if and only if there are nonzero vectors $v_0$ of
size $(k_0+1) \times 1$ and $v_1$ of size  $(k_1+1) \times 1 $ both with constant
coefficients such that 
\begin{equation}
  \label{eq:00}
Av_0=v_1 h;     
\end{equation}
the  tensor  $ \Hh=v_0 \otimes v_1$ is  called  an unstable (or jumping) hyperplane for the matrix $A$.  

For $p\geq3$ there are at least two ways to define a jumping
hyperplane. We will call them weak and strong jumping hyperplanes. 
\begin{defin}
 \label{def:wjumping}
 $\Hh=v_0 \otimes v_j \otimes h \in V_0 \otimes V_j \otimes \widehat{V}^j$
    (where $\widehat{V}^j=V_1 \otimes \dots \otimes \widehat{V_j} \otimes
    \dots \otimes V_p$ ) is a {\bf $(j)$-weak jumping hyperplane} for $A$ if $\exists \ v_0,w_1, \dots , w_{k_0}$
  basis of $V_0$ such that
  \begin{equation}
  \label{eq:wjump}
   A=v_0 \otimes v_j \otimes h +\sum_{i=1}^{k_0}w_i \otimes \dots
  \end{equation}
  where $h\in \widehat{V}^j$ generate an hyperplane for $\PP^{k_1} \times \dots \times
  \widehat{\PP^{k_j}} \times \dots \times \PP^{k_p} \subset \PP(V_1
  \otimes \dots \otimes \widehat{V_j} \otimes \dots \otimes V_p)$
  (that, by abuse of notations, we call also $h$.)

\end{defin}
\begin{remark}
  The expression ({\ref{eq:wjump}}) means, as in the case $p=2$, that
  $H^0(Ker{f_A}^{(j)}_{\vert h})\!\!\neq\!\!0.$ (i.e. , by definition,
  $h$ is a jumping hyperplane for the bundle ${S_A}^{(j)}$).
\end{remark}

 If $\Hh= v_0 \otimes v_j $ is a $(j)$-weak jumping hyperplane for $A$ then the map:
 \begin{center}
\begin{equation*}
     \begin{aligned}
   V_0 \otimes \dots \otimes V_p & \to (V_0 /_{<v_0>}) \otimes \dots \otimes (V_j/_{<v_j>}) \otimes \dots  \otimes V_p \\
   A & \mapsto {A'}_{j} \\
 \end{aligned}
 \end{equation*}
\end{center}
gives an elementary transformation \cite{M1}.

\begin{remark}
${A'}_j$ is again of boundary format. In particular, after a basis
has been chosen, ${A'}_j$ is obtained  by deleting two
directions in $A$ .  
\end{remark}
\begin{prop}
If ${A'}_j$ is defined as above  $$ Det A \neq 0  \Rightarrow Det {A'}_j \neq 0 $$
\end{prop}
\begin{proof}
  If $X:=\PP^{k_1} \times \dots \times \widehat{\PP^{k_j}} \times
  \dots \times \PP^{k_p}$ and $h$ is the hyperplane defined in
  (\ref{def:wjumping}) associated to $\Hh$, the map ${S_A}^{(j)} \to
  \OO_h $ induced by a non zero section of ${S_A}^{(j)}$ is surjective
  (the same proof of \cite{val} prop.2.1 works).

\noindent Since $\codim h=1$, then its kernel
  ${S'}^{(j)}$ is locally free sheaf \cite{ser} of rank $k_0-k_j-1$ on
  $X$ and it is the Steiner bundle associated to the matrix ${A}^{(j)}$ as the
   snake-lemma applied to the following exact diagram shows
   \begin{small}
     \begin{equation*}
    \label{eq:diagram}
      \begin{CD}
     @. @.  @. 0\\
    @.  @. @. @VVV \\      
    @. @.  @.  {S'}^{(j)}\\
    @.  @. @. @VVV \\      
  0 @>>>\OO_X(-1, \dots, -1) \otimes {V_j}^{\ch}
       @>{f_A}^{(j)}>>\OO_X \otimes V_0^{\ch} @>>>{{S_A}}^{(j)} @>>>0 \\
      @. @VVV @VVV @VVV   \\
    0 @>>> \OO_X(-1, \dots, -1) @>>>\OO_X @>>> \OO_h @>>>0 \\
       @. @VVV @VVV @VVV   \\
     @. 0  @. 0 @. 0
 \end{CD}
    \end{equation*}
\end{small}
i.e. , ${S'}^{(j)}= {S_{A'_{j}}}^{(j)}$
and by proposition (\ref{prop:det}) the result follows.
\end{proof}
\begin{remark}
  If $W(S_{A}^{(j)})$ is the set of jumping hyperplanes of
  the bundle $S_{A}^{(j)}$ ,
  then the exact sequence (dual to the last column of the above
  diagram)
  \[0 \to {S_{A}^{(j)}}^\ch \to{S_{A'_{j}}^{(j)}}^\ch \to
  \OO_X(1,\dots,1)\to0 \]
shows that $W({S_{A}^{(j)}}) \subset W(S_{A'_{j}}^{(j)}) \cup \{h\}$ 
\end{remark}
%%%remark sull'azione e le trasf. elementari
%%%\begin{remark}
%%%  An element $ g \in Stab(A) $ preserves $h$ and it
%%%  induces $\bar{g} \in SL(V_0/_{<g(v_0)>}) \times SL(V_1) \times \dots
%%%  \times SL(V_j/_{<g(v_j)>})\times \dots \times SL(V_p)$
%%%  such that $g \cdot A$ projects to $\bar{g} \cdot {A'}_j$ and the elementary
%%%  transformation behaves well with respect to the action of $g$.
%%%\end{remark}
\begin{defin}
\label{def:sjumping}
  $\Hh=v_0 \otimes v_1 \otimes \dots \otimes v_p$ is a {\bf strong
    jumping hyperplane} for $A$ if $ \exists \ v_0,w_1, \dots ,
  w_{k_0}$ basis of $V_0$ such that
\begin{equation*}
A = v_0 \otimes v_1 \otimes \dots \otimes v_p + \sum_{i=1}^{k_0} w_i \otimes \dots
\end{equation*}
\end{defin}
\begin{remark}
  If $\Hh$ is a strong jumping hyperplane then $\Hh$ defines a $(j)$-weak
  jumping hyperplane for all $j=1,\dots,p$; in particular for a
  strong jumping hyperplane there are many elementary transformations.
\end{remark}
\begin{remark}
  For $p=2$ the notations of strong jumping hyperplane and of weak jumping
  hyperplane coincide with each other (see \cite{ao}).
\end{remark}

\begin{exm} ({\it the identity})
Fixed a basis $e_0^{(j)}, \dots , e_{k_j}^{(j)} $ in $V_j$ for
all $j$, the identity matrix is represented by
\begin{equation*}
I:= \sum_{\substack{
i_0=i_1+\dots +i_p \\
 0 \leq i_j \leq k_j}} e_{i_0} ^{(0)} \otimes \dots \otimes e_{i_p}^{(p)}.
\end{equation*}

Let $t_0, \dots, t_{k_0}$ be any distinct complex numbers. Let $w$ be
the $(k_0+1)\times (k_0+1)$ Vandermonde matrix whose $(i,j)$ entry is
$t_j^{(i-1)}$, so acting with $w$ over $V_0$, we have:
\begin{equation*}
 e_j^ {(0)}=\sum _{s=0}^{k_0} \bar{e}_s^{(0)}t_s^j
\end{equation*}
Then substituting

\begin{equation*}
\begin{split}
  I &= \sum_{\substack{
i_0=i_1+\dots +i_p \\
  s=0,\dots k_0
}}  \bar{e}_{i_0} ^{(0)}t_s^{i_0} \otimes e_{i_1}^{(1)} \otimes \dots \otimes e_{i_p}^{(p)} \\
&=\sum_{s=0}^{k_0}\bar{e}_s^0 \otimes (\sum_{i_1=0}^{k_1}
e_{i_1}^{(1)}t_s^{i_1}) \otimes \dots \otimes (\sum_{i_p=0}^{k_p}
e_{i_p}^{(p)} t_s^{i_p})
\end{split}
 \end{equation*}
 
Thus, since $t_i$ have no restrictions, $I$ has infinitely many strong
 jumping hyperplane.\\
 We call {\it Schwarzenberger bundle} the vector bundle associated to $I$
 (in fact in the case $p=2$ it is exactly the same introduced by
 Schwarzenberger in \cite{s})(see also (\cite{ao})
\end{exm}
\begin{prop}
\label{prop:identity}
Let $A$ be a boundary format matrix with $Det A \neq 0$.
If $A$ has $N \geq k_0+3 $ strong jumping hyperplanes then it is an identity.
\end{prop}
\begin{proof}
  In the case $p=2$ the statement is proved in \cite{ao} (theorem
  5.13) or in \cite{val} (theorem 3.1).
  Chosen $V_0$ and other two vector spaces among $V_1, \dots, V_p$
  (say $V_1$ and $V_2$), one may perform several elementary
  transformations with $V_0$ and all the others so that we get $A' \in
  {V'}_0 \otimes V_1 \otimes V_2 $ boundary format matrix with $Det A'
  \neq 0$ and $N' \geq {k'}_0 +3 $ strong jumping hyperplanes, then 
  $A'$ is an identity.\\
  As in the above example,
  one can change the hyperplane giving the elementary transformation,
  so that for all $N$ strong jumping hyperplanes we get $t_1, \dots,
  t_N$ distinct complex numbers and corresponding suitable basis of
  $V_1$ and $V_2$ :
\begin{equation*}
\begin{aligned}
  {{\bar e}_0}^{(1)}& \dots {{\bar e}_{k_1}}^{(1)}\\
  {{\bar e}_0}^{(2)}& \dots {{\bar e}_{k_2}}^{(2)}\\
\end{aligned}
\end{equation*}
such that the hyperplanes are given by
\begin{equation*}
\sum_{i=0}^{k_1}{\bar e}_i^{(1)}t_j^i \quad \text{and} \quad
\sum_{i=0}^{k_2}{\bar e}_i^{(2)}t_j^i \quad \text{for} j=1, \dots N
\end{equation*}
Now, changing  $V_1$ and $V_2$ with the pairs $V_1, V_j$ ($j=1, \dots p$ ) we get
\begin{equation*}
A:=\sum_{s=0}^{k_0}\bar{e}_s^0 \otimes (\sum_{i_1=0}^{k_1}
e_{i_1}^{(1)}t_s^{i_1})  \otimes \dots \otimes (\sum_{i_p=0}^{k_p}
 e_{i_p}^{(p)} t_s^{i_p})
\end{equation*}
showing that $A$ is an identity.
\end{proof}

\begin{prop}
Two nondegenerate boundary format matrices having in common
$k_0+2$
distinct strong jumping hyperplanes determine isomorphic Steiner
bundles for every $j$.
\end{prop}
\begin{proof}
 In the case $p=2$ the statement is proved in \cite{ao} (theorem 5.3).
  Chosen $V_0$ and other two vector spaces among $V_1, \dots, V_p$
  (say $V_1$ and $V_2$), one may perform several elementary
  transformations with $V_0$ and all the others so that we get $A' \in
  {V'}_0 \otimes V_1 \otimes V_2 $ boundary format matrix with $Det A'
  \neq 0$ and $N' = {k'}_0 +2 $ strong jumping hyperplanes, then 
  $S_{A'}^{(j)}$ is uniquely determined for every $j$. Now, changing $V_1$ and $V_2$
 with the pairs $V_1$ and $V_j$ $(j=2, \dots, p)$ we detect  all the 
$3$-dimensional submatrices of $A$ which give bundles  uniquely determined, 
so also $S_A^{(j)}$ is uniquely determined for every $j$.
\end{proof}

\begin{remark}
  In the case $p=2$ we know that $k_0+2$ jumping hyperplanes give an
  existence condition for the bundles $S_A^{(j)}$ (they are  logarithmic
  bundles, see \cite{ao}) but in the case $p \geq 3$ there is not an analog
  existence result.(The previous proposition gives only the uniqueness)
\end{remark}
The following is a classical result (see for instance \cite{ha}
prop.9.4 page 102, or \cite{dk} theorem 6.8)
\begin{prop} \label{represented}
All nondegenerate matrices of type $2\times k\times (k+1)$ are
$GL(2)\times GL(k)\times GL(k+1)$ equivalent, or equivalently 
every surjective morphism of vector bundles on $\PP^1$
\[\OO_{\PP^1}^{k+1}\to\OO_{\PP^1}(1)^k\]
is represented by an identity matrix.
\end{prop}

We recall now the following
\begin{prop}\cite{ao}
 \label{def:identity}
Let 
$A \in V_0 \otimes  \dots \otimes V_p$ A be a $(p+1)$-dimensional
matrix of boundary format the following conditions are equivalent:
\begin{itemize}
       \item[(i)] $A$ is an identity;
      \item[(ii)] there exist a vector space $U$ of dimension 2 and
       isomorphisms $V_j \simeq S^{k_j}U$ such that $A$ belongs to the
       unique one dimensional $SL(U)$-invariant subspace of $S^{k_0}U
       \otimes \dots \otimes S^{k_p}U$.
     \end{itemize} 
\end{prop}
The equivalence between (i) and (ii) follows easily from the following remark:
 the matrix  $A$ satisfies the condition (ii) if and only if
it corresponds to the natural multiplication map
$S^{k_1}U\otimes\ldots\otimes S^{k_p}U\to S^{k_0}U$ (after a suitable
isomorphism $U\simeq U^\ch$ has been fixed). We notice that by the
 Clebsch-Gordan decomposition of the tensor product there is a unique
 $SL(U)$-invariant map as above.
 
\begin{remark}
  If $A$ is not an identity, an element $g \in Stab(A) $ preserves a
  $(j)$-weak jumping hyperplane $h$ and it
  induces $\bar{g} \in SL(V_0/_{<g(v_0)>}) \times SL(V_1) \times \dots
  \times SL(V_j/_{<g(v_j)>})\times \dots \times SL(V_p)$
  such that $g \cdot A$ projects to $\bar{g} \cdot {A'}_j$ and the elementary
  transformation behaves well with respect to the action of $g$.
\end{remark}

 For every integer $j$, let  $\Djs(A)$ be  the {\it locus of
 $(j)$-strong directions  of $A$} defined as 
\begin{equation*}
  \begin{split}
\{ <v_j>\in \PP(V_j^\ch) \  \vert  \ & \forall i\neq j \ \exists \ v_i \in V_i 
   \ \text{such that}\\
\ &  v_0 \otimes \dots
   \otimes  v_p \ \text{is a strong jumping hyperplane for} A \}
\end{split}
\end{equation*}

We racall that (see for details  \cite{ao}) for boundary format
matrices  the following conditions
are equivalent
\begin{enumerate}
\item $A\in \Vop$ is diagonal
\item $\CC^* \subset Stab (A)$
\item there exist a vector space $U$ of dimension 2, a subgroup
  $\CC^* \subset SL(U)$ and isomorphisms $V_j \simeq S^{k_j}U$ such
  that $A$ is a fixed point of the induced action of $\CC^*$.
\end{enumerate}
Then, the same proofs of  corollaries 6.9 -
6.10 and lemmas 6.12-6.13 of \cite{ao} work  also in
the $(p+1)$-dimensional case, by replacing $V$ by $V_j$ and $W(S)$ by
$\Djs(A)$. More precisely we have:
\begin{cor}
  Let $A$ be a boundary format  nondegenerate matrix. If $\CC^*
  \subset Stab(A)$
  then  for every $j$ the $\CC^*$-action on $V_j$ has exactly $k_j+1$
  fixed points whose weights are proportional to  $-k_j, -k_j+2,\dots,
  k_j-2, k_j$.  
\end{cor}
\begin{remark}
\label{actionplus} More in general,
 the $\CC^*$-action on $V$ (where $V$ is a $n+1$-dimensional
  vector space) has exactly $n+1$
  fixed points whose weights are proportional to  $-n, -n+2,\dots,
  n-2, n$ if and only if  there exist a vector space $U$ of dimension
  2 such that $\CC^* \subset SL(U)$ and  $V \simeq S^{n}U$.
\end{remark}
\begin{cor}
  Let $A$ be a boundary format  nondegenerate matrix such that $\CC^* \subset
  Stab(A)$. Then either $A$ is an identity or 
$\Djs(A)$ has only two closed points, namely   the two fixed points 
of the dual $\CC^*$-action on $\PP(V_j^\ch)$  having minimum and maximum weights.
\end{cor}

\begin{lemma}
  Let $U$ be a $2$-dimensional vector space, and $\forall j$ \ 
  $C_j\simeq\PP(U)\to\PP(S^{k_j}U)$ be the $SL(U)$-equivariant
  embedding (whose image is a rational normal curve).  Let
  $\CC^*\subset SL(U)$ act on $\PP(S^{k_j}U)$.  We label the $k_j+1$
  fixed points $P_i, \ i=-k_j+2n,\ n = 0,\dots, k_j$ of the
  $\CC^*$-action with an index proportional to its weight.  Then
  $P_{-k_j}$, $P_{k_j}$ lie on $C_j$ and $P_{-k_j+2n}=T^nP_{-k_j}\cap
  T^{k_j-n}P_{k_j}$, where $T^n$ denotes the $n$-dimensional
  osculating space to $C_j$.
\end{lemma}

\begin{lemma}
  \label{twodiff}
  Let $A$ be a boundary format nondegenerate matrix.
  If there are two different one-parameter subgroups
  $\lambda_1, \lambda_2: \CC^* \to Stab(A)$ then $A$ is an identity.
\end{lemma}

\vspace{0.5cm}
{\bf Proof of theorem \ref{th:main}}
\begin{proof}
  We proceed by induction on $k_0$.

  \noindent  If $k_0=2$ the theorem is true by
  proposition ~\ref{represented}.

\noindent When $Stab(A)^0$ contains only the
  identity  the result is trivial hence we may suppose that 
   $\dim Stab(A)^0 \geq 1$ then, according to  (\cite{ao}, theorem $2.4$)
  the matrix $A$
  is triangulable and  there exists at least one 
   strong jumping hyperplane $\Hh=v_0\otimes \dots \otimes v_p$.\\
  We may  also suppose that the number of jumping hyperplanes is
  finite otherwise A is an identity (proposition \ref{prop:identity}),
  hence $\Hh$ is $Stab(A)^0$- invariant.  Let $A_1'$ be the image of $A$
  by the elementary transformation associated to the $(1)$-weak jumping
  hyperplane  defined by $\Hh$ (we choose $j=1$ to
  have simpler notations). The matrix $A_1'$ belongs to $V_0' \otimes
  V_1' \otimes V_{2} \otimes \dots \otimes V_{p}$ where $V_0'=V_0
  /_{<v_0>}$ and $V_1'=V_1/{<v_1>}$, it is nondegenerate and of
  boundary format then, by induction, there exists a $2$-dimensional
  vector space $U$ such that
\[ V_0'  \simeq S^{k_0-1}(U),\quad  V_1'  \simeq S^{k_1-1}(U) \quad \text{and} 
\quad V_i=S^{k_i}(U) \quad \text{for all} \quad i \geq 2 \] and
$Stab(A_1')^0 \subseteq SL(U)$ (by using essentially the same
argument we could work in $GL(V_0) \times \dots \times GL(V_p)$).

\noindent  Since $A_1'$ is obtained  from the matrix $A$ after the choice of
two directions, any element which stabilizes $A$ also stabilizes
$A_1'$, so $Stab(A)^0 \subseteq Stab(A_1')^0$.
Hence $Stab(A)^0 \subseteq SL(U)$ and $SL(U)$ acts on $V_i$ according
to $V_i \simeq S^{k_i}U$ for $i\geq 2$, by the inductive hypothesis. 

\noindent Now, we claim that the action of $SL(U)$ can be lifted
to the whole $\Vop$.

\noindent Indeed, the above considered elementary
transformation gives the decomposition $V_0=V_0' \oplus \CC$ and $V_1
=V_1' \oplus \CC$. 

\noindent
If $\phi: \CC^* \to GL(V_i')$ is the natural action of $\CC^*\subset
SL(U)$ on $V_i'=S^{k_i-1}U$ (for $i=0,1$) with $k_i$ fixed points
having weights $-k_i+1, -k_i+3, \dots, k_i-1$,  we can construct an
action  $\psi: \CC^* \to  GL(V_i' \oplus \CC)$ on $V_i$ defined by
\[ t \mapsto 
\begin{pmatrix}
  t^{-1} \phi(t)& 0 \\
   0 & t^{k_i}
\end{pmatrix}
\]
 having $k_i+1$ fixed points with weights $-k_i,
-k_i+2, \dots, k_i$.
hence, by remark \ref{actionplus}, the statement  follows.

\noindent In the case $Stab(A)^0=SL(2)$, the action of $SL(U)$
satisfies definition \ref{def:identity}, proving that $A$ is an
identity.  
\noindent Now, as in \cite{ao}, consider the Levi
decomposition $Stab(A)^0=M \cdot R$ where R is the radical and M is
maximal semisimple. If $A$ is not an identity (i.e. , $Stab^0(A) \neq
SL(2)$) then $M=0$ and $Stab(A)^0$ is solvable hence by the Lie theorem it
is contained (after a convenient basis  has been chosen) in the subgroup of
upper triangular matrices $T=\left\{
  \begin{pmatrix}
    a&b\\ 0&\frac{1}{a}\\
  \end{pmatrix}
\vert \ a\in\CC^*, b\in\CC\right\}$. If there is a subgroup $\CC^*$
properly contained in $Stab(A)^0$ then there is a conjugate of $\CC^*$ different from
itself and this is a contradiction by the lemma \ref{twodiff}.
If 
$Stab(A)^0$ does not contain proper subgroups $\CC^*$ then it is isomorphic to
$\CC\simeq\left\{
\begin{pmatrix}
    1&b\\ 0&1\\
  \end{pmatrix}
\vert  \ b\in\CC\right\}$.
\end{proof}
\begin{remark}
  Throughout this paper we work only on nondegenerate matrices.
  Indeed, in the proofs we apply the induction strategy (hence the
  results of \cite{ao}) and the correspondence between matrices and
  vector bundles described in proposition \ref{prop:det}.
 
 \noindent 
 The characterization of the stabilizer of degenerate matrices is
 still an open problem. \\
 Another interesting problem is the study of the stabilizer of general
 multidimensional matrices (and not necessarily of boundary format).

\end{remark}


\begin{thebibliography}{Bar77b}
                                                            
\bibitem[AO]{ao} V.~Ancona and G.~Ottaviani, \emph{{Unstable
      hyperplanes for Steiner bundles and multidimensional matrices}},
  Advances in Geometry, 1 (2001), 165-192
 

\bibitem[D]{cd} C.~Dionisi, \newblock {\em Multidimensional
  matrices and minimal resolutions of vector bundles} \newblock Ph.D. thesis, Dip. Matematica "R.Caccioppoli", Universit\`a di Napoli,
  2000.

\bibitem[DO] {do} C.~Dionisi and G.~Ottaviani, \emph{The Binet-Cauchy
  theorem for the hyperdeterminant of boundary format multidimensional
  matrices}, Journal of Algebra, 259, (2003),
  87-94

\bibitem[DK] {dk}
I.~Dolgacev and M.~M.~Kapranov, \emph{{Arrangement of hyperplanes and
  vector bundles on $\PP^n$}}, Duke Math. J. 71 (1993), 633-664

\bibitem[GKZ]{gkz}
I.~M.~Gelfand, M.~M.~Kapranov and A.~V.~Zelevinski, \emph{{Discriminants,
resultants and multidimensional determinants}}, Birkhaeuser, Boston 1994

\bibitem[GKZ1]{gkz1}
I.~M.~Gelfand, M.~M.~Kapranov and A.~V.~Zelevinski,
\emph{{Hyperdeterminants}}, Adv. in Math. (1992), no.~96, 226--263

\bibitem[Ha]{ha}
J.~Harris, \emph{{Algebraic geometry, a first course}}, Graduate Texts in
  Mathematics, no. 133, Springer-Verlag, New York-Heidelberg-Berlin,
  1992.

\bibitem[Mar82]{M1}
M.~Maruyama, \emph{{Elementary tranformations in the theory of algebraic vector
  bundles}}, Lect. Notes Math. (1982), no.~961, 241--266.

\bibitem[Par]{parf}
P.~G.~Parfenov, \emph{{Orbits and their closures in the spaces
    $\CC^{k_1} \otimes \dots \otimes \CC^{k_r}$}}, Sbornik
Math. (2001), \textbf{192}, 89--112 (English transl.)

\bibitem[Sch]{s}
R.L.E. Schwarzenberger, \emph{{Vector bundles on the projective
    plane}}, Proc.London Math.Soc. \textbf{3} (1961), no.~11, 623--640.

\bibitem[Ser65]{ser}
J.P. Serre, \emph{{Alg\`ebre locale, multiplicit\'es}}, Lecture Notes in
  Mathematics, no.~11, Springer-Verlag, 1965.


\bibitem[WZ]{wz}
J.~Weyman and A.~V.~Zelevinsky, \emph{{Singularities of
    hyperdeterminants}}, Ann. Inst. Fourier 46 (1996), 591-644.

\bibitem[V1]{val1} J.~Vall\`es,  \emph{Fibr\'es de Schwarzenberger et coniques
  de droites sauteuses}, Bull.~Soc.~Math.~France, \textbf{128}, (2000), 433-449.

\bibitem[V2]{val} J.~Vall\`es, \emph{Nombre maximal d'hyperplans instables pour
  un fibr\'e de Steiner}, Math. Zeitschrift 233, 507-514 (2000).
\end{thebibliography}
\end{document}